\documentclass[12pt]{article}
\usepackage[T2A]{fontenc}
\usepackage[english]{babel}
\usepackage{amsfonts,amssymb,amsmath,enumerate,amsthm}
\usepackage{epsfig}
\usepackage{float}

\newcommand{\RR}{\mathbb R}
\newcommand{\CC}{\mathbb C}
\newcommand{\NN}{{\mathbb N}}

\newtheorem{theorem}{Theorem}
\newtheorem{lemma}{Lemma}
\newtheorem{remark}{Remark}

\newtheorem{corollary}{Corollary}

\newcommand{\beq}{\begin{equation}}
\newcommand{\eeq}{\end{equation}}
\newcommand{\ba}{\begin{array}}
\newcommand{\ea}{\end{array}}
\newcommand{\bea}{\begin{eqnarray}}
\newcommand{\eea}{\end{eqnarray}}

\providecommand{\keywords}[1]
{
  \small	
  \textbf{Keywords: } #1
}
\providecommand{\AMS}[1]
{
  \small	
  \textbf{AMS Subject classification: } #1
}

\usepackage{graphicx}
\usepackage[usenames,dvipsnames]{color}
\usepackage{transparent}

\DeclareMathAlphabet{\mathpzc}{OT1}{pzc}{m}{it}

\usepackage{ulem}
\usepackage{cancel}

\begin{document}
\begin{center}

{\bf A holographic uniqueness theorem for the two-dimensional Helmholtz equation}
\vskip 10pt
{\it A.V. Nair,  \it R.G. Novikov }

\vskip 10pt

\end{center}

\begin{abstract}
We consider a plane wave, a radiation solution, and the sum of these
solutions (total solution) for the Helmholtz equation in an exterior region in $\RR^2$.
We consider a straight line in this region, such that the direction of propagation of the  plane wave 
is not parallel to this line.
We show that the radiation solution in the exterior region is uniquely determined by the intensity of 
the total solution on an interval of this line.
In particular, this result solves one of the old mathematical questions of holography in its two-dimensional setting. Our proofs also contribute to the theory of the Karp expansion of radiation solutions in two dimensions.
\end{abstract} 
\AMS{35J05, 35P25, 35R30}
\section{Introduction}
We consider the two-dimensional Helmholtz equation
\begin{equation}
\label{eq:1.1}
\Delta\psi(x) + \kappa^2\psi(x) = 0, \ \  x \in {\cal U},\ \  \kappa>0, 
\end{equation}
where $\Delta$ is the Laplacian in $x$, and $\cal U$ is a region (open connected set) in $\RR^2$
consisting of all points outside a closed bounded regular curve $S$ (as in \cite{K}). 
This equation particularly arises in electrodynamics, acoustics, and quantum mechanics.
For equation (\ref{eq:1.1}) we consider solutions $\psi_0$ and $\psi_1$ such that:
\begin{equation}
\label{eq:1.2}
\psi_0 = e^{ikx}, \ \ k \in \mathbb{R}^2, \ \ |k| = \kappa,
\end{equation}
$\psi_1$ is of class $C^2$ and satisfies the Sommerfeld's radiation condition
\begin{equation}
\label{eq:1.3}
\sqrt{|x|}\bigl(\frac{\partial}{\partial|x|} - i\kappa \bigr)\psi_1(x) \to 0 \ \ as \ \ |x| \to +\infty,
\end{equation}
uniformly in $x/|x|$. We say that $\psi_0$ is the plane wave solution and $\psi_1$ is a radiation solution.

Let
\begin{equation}
L = L_{x_0, \theta}=\{x\in\RR^2:\ \ x=x(s)=x_1+s\theta,\ \ -\infty<s <+\infty  \}, \nonumber
\end{equation}
\begin{equation}  
 L^+=L^+_{x_1, \theta}=\{x\in\RR^2:x=x(s)=x_1+s\theta,\ \ 0<s <+\infty  \}, \nonumber
\end{equation}
\begin{equation}
\label{eq:1.4}
 L^-=L^-_{x_2, \theta}=\{x\in\RR^2:\ \ x=x(s)=x_2-s\theta,\ \ 0<s <+\infty  \},
\end{equation}
where $x_0,x_1,x_2 \in \RR^2$, $\theta \in \mathbb S^1$ and $x_1,x_2 \in L$. \\ Thus, $L=L_{x_0, \theta}$ is the oriented straight line in  $\RR^2$  that passes through $x_0$ and has direction $\theta$.  $L^+=L^+_{x_1, \theta}$ is the ray in  $\RR^2$  that starts at $x_0 \in\RR^2$ and has direction  $\theta$, and $L^- = L^-_{x_2, \theta}$ a ray that lies on the extension of $L^+$ in the $(-\theta)$ direction.

In the present work we show that, for a fixed plane wave solution $\psi_ 0$, 
any complex-valued radiation solution $\psi_ 1$  on $L^+ \cup L^-$ is uniquely determined by the intensity $|\psi|^2$  of the total solution $\psi= \psi_ 0+\psi_ 1$ on arbitrary intervals $\Lambda^+ \cup \Lambda^-$of $L^+$ and $L^-$ respectively,
under the assumptions that $L^+=L^+_{x_1, \theta}\subset ~{\cal U}$, $L^-=L^-_{x_2, \theta}\subset ~{\cal U}$ and $\pm\theta \neq k/|k|$. Here,  ${\cal U}$ is the region in (\ref{eq:1.1}) and $k$ is the vector in (\ref{eq:1.2}).  
This result is given as Theorem~1 in Section 2. 

As a corollary, we also obtain that, \\
for any straight line $L\subset  {\cal U}$ such that $\pm\theta \neq k/|k|$, any complex-valued  $\psi_ 1$  on $L$ is uniquely determined  by  the  intensity $|\psi_ 0+\psi_ 1|^2$
on an arbitrary interval $\Lambda$ of $L$,  for  fixed   $\psi_ 0$; see Theorem 2 in Section 2,

Our studies are motivated by problems of holography and phaseless inverse scattering.
These research areas go back, in particular,  to  \cite{D}, \cite{E}, \cite{EW}.
In connection with recent or relatively recent mathematical results obtained in these directions,
see, for example,
 \cite{HN}, \cite{JL}, \cite{MV} - \cite{GR}, \cite{E}, \cite{EW},
and references therein.

The present work continues \cite{N} where holography uniqueness results were proved for the Helmholtz equation in the three-dimensional case.

Theorems 1, 2  \ mentioned above and presented in detail in Section 2 solve one of the old mathematical questions arising in holography in its two-dimensional setting. 

In the present work we use the Karp expansion (\ref{eq:3.2}) (recalled in Section 3) instead of the Atkinson-Wilcox expansion used in \cite{N} for radiation solutions of the Helmholtz equation. Although the Karp expansion can be considered as a two-dimensional analogue of the three-dimensional Atkinson-Wilcox expansion used in \cite{N}, its properties are somewhat different, more complicated, and less studied. In this connection the present work also contributes to the theory of Karp expansion of the radiation solutions $\psi_1$ of the two-dimensional Helmholtz equation (\ref{eq:1.1}); see Lemmas 2, 3, and Theorem 3 in Section 4, and proof of Lemma 1 in Section 7. 

The main results of this work are presented in more detail and proved in Sections 2, 4, 5, 6, 7.
In our proofs, we proceed from the results recalled in Section~3.

\section{Main results}

Our key result is as follows.

\begin{theorem}\label{thm:1}
Let $\psi_0=e^{ikx}$ and  $\psi_1$  be solutions of equation (\ref{eq:1.1}) as in formulas (\ref{eq:1.2}) and (\ref{eq:1.3}).
Let $L$, $L^+$ and $L^-$ be as given in (\ref{eq:1.4}) such that  $L^+=L^+_{x_1, \theta}\subset {\cal U}$ and $L^-=L^-_{x_2, \theta}\subset {\cal U}$, $x_1,x_2 \in L$ , $\theta \neq k/|k|$, and $-\theta \neq k/|k|$, where ${\cal U}$ is the region in  (\ref{eq:1.1}). Then $\psi_1$ on $L^+$ $\cup$ $L^-$ is  uniquely determined by the  intensity  $|\psi_0+\psi_1|^2$ on   $\Lambda^+\cup \Lambda^-$,  for fixed $k$,
where $\Lambda^+$ and $\Lambda^-$ are arbitrary non-empty intervals of $L^+$ and $L^-$ respectively.
\end{theorem}

As a corollary, we also get the following result.

\begin{theorem}\label{thm:2}
Let $\psi_0=e^{ikx}$ and  $\psi_1$  be solutions of equation (\ref{eq:1.1}) as in formulas (\ref{eq:1.2}) and (\ref{eq:1.3}).
Let $L$ be a straight line in $\RR^2$ such that $L\subset {\cal U}$ as described in (\ref{eq:1.4}). If $\theta \neq k/|k|$ and $-\theta \neq k/|k|$, then $\psi_1$ on $L$ is  uniquely determined by the  intensity  $|\psi_0+\psi_1|^2$ on  $\Lambda$,  for fixed $k$,
where  $\Lambda$  is an arbitrary non-empty interval of $L$.
\end{theorem}
Theorem 1 is proved in Sections 4, 5, and 6. This proof uses techniques developed in \cite{K} and \cite{N}. Note that $\psi_0$, $\psi_1$ are real-analytic on ${\cal U}$, and, therefore, on $L^+$, $L^-$ in Theorem~1 and on $L$ in Theorem 2. Therefore,  the function $|\psi|^2=|\psi_0+\psi_1|^2= (\psi_0+\psi_1)(\overline \psi_0+\overline \psi_1)$ is  real-analytic on $L^+$, $L^-$ in Theorem 1 and on $L$ in Theorem 2. Given this analyticity,  Theorem~1 reduces to the case when  $\Lambda^+ = L^+$, $\Lambda^- = L^-$ and Theorem 2 reduces to the case when $\Lambda=L$. 

Theorem 2 is proved as follows (for example). We assume that $\Lambda=L$. Then we simply consider $L^+ \subset L$ and $L^-\subset L$ such that $L^+ \cap L^- \neq \emptyset$. Then we may apply theorem 1. 

\begin{corollary}
Under the assumptions of Theorem 2,  the  intensity  $|\psi_0+\psi_1|^2$ on  $L$, for fixed $k$, uniquely determines $\psi_1$  
in the entire region ${\cal U}$.  

\end{corollary}
Corollary 1 follows from  Theorem 2, formula (\ref{eq:3.9}), and analyticity of  $\psi_1$ in ${\cal U}$.

For the three-dimensional Helmholtz equation, prototypes of Theorems 1,2 \ and Corollary 1 \ were obtained in \cite{N}. In the present work, we proceed from Karp's results on two-dimensional radiation solutions (in \cite{K}) for equation (\ref{eq:1.1}) and the approach of \cite{N}.

Note that the problem in two dimensions is more difficult than in three dimensions. The reason is that the Karp expansion (\ref{eq:3.2}) for the radiation solution in two dimensions is considerably more complicated than the Atkinson-Wilson expansion (in \cite{A}, \cite{W}) for the three-dimensional case.

In particular in Theorem 1 \ we need to consider two rays $L^+$ and $L^-$ in place of a single ray $L^+$ as was in \cite{N}. Both these rays need to have different directions than the propagation direction of the plane wave $\psi_0 = e^{ikx}$. In Theorem 2 \ as well we need both the directions of the line $L$ to be different from the propagation direction of the plane wave $\psi_0 = e^{ikx}$.

In the present work, we strongly use the uniqueness of Karp's expansion (\ref{eq:3.2}), which is formulated as Lemma 1 \ in the next section. For completeness of presentation, we prove this lemma, as we did not find a proof in the literature.

In the present work we also contribute to the theory of the Karp expansion for radiation solutions of the two-dimensional Helmholtz equation (\ref{eq:1.1}) by Lemmas 2, \ 3, \ and Theorem 3 \ given in Section 4. In particular, in Lemma 2 \  we write explicitly recursive relations for coefficients of Karp expansion, in Lemma 3 \ we establish important symmetries of these coefficients and in Theorem 3, we present a recursive method for finding all coefficients in Karp expansion from coefficients arising in standard scattering theory.

Note also that Theorems 1,\ 2,\  and Corollary 1 \ admit straightforward applications to phaseless inverse scattering in two dimensions. 
\section{Preliminaries}
Let $(r,\phi)$ denote the polar coordinates of a point $x$ and let $\theta = \frac{x}{|x|}$ be its direction. So $-\theta$ will be the direction of $-x$ whose polar coordinates will be $(r, \phi + \pi)$. \\
We use the Laplacian in polar coordinates:
\begin{equation}
\label{eq:3.8}  
\Delta = \frac{\partial^2}{\partial r^2} + \frac{\partial}{r\partial r} + \frac{\partial^2}{r^2\partial\phi^2} \ .
\end{equation}
Let
\beq
\label{eq:3.1}
B_{\rho}=\{x\in\RR^2 \ | \ |x|<\rho\},\ \ \rho>0. 
\eeq
Suppose that  $\psi_1$ is a radiation solution of  equation (\ref{eq:1.1}), and  $\RR^2\setminus B_{\rho}  \subset {\cal U}$. \\
Due to Theorems I, II \  of \cite{K},  we have that: 
\begin{equation}
\label{eq:3.2}  
\psi_1(x)=H_0(\kappa |x|)\sum\limits_{j=0}^{\infty}\frac{F_j(\phi)}{|x|^{j}} + H_1(\kappa |x|)\sum\limits_{j=0}^{\infty}\frac{G_j(\phi)}{|x|^{j}},
\end{equation}
\begin{center}
{$\mbox{\rm for}\ \ x\in \RR^2\setminus B_\rho,$}
\end{center}
where $H_0$ and $H_1$ are the Hankel functions of the first kind of order zero and one respectively. The series converges absolutely and uniformly for $|x| \geqq \rho_1 > \rho$ (for any such $\rho_1$).
\begin{lemma}
Suppose that: 
\begin{equation}
\label{eq:3.8}
H_0(r)\sum\limits_{j=0}^{\infty}\frac{A_j}{r^{j}} + H_1(r)\sum\limits_{j=0}^{\infty}\frac{B_j}{r^{j}} = 0,
\end{equation}
for $r \geq r_1 > 0$, $A_j, B_j \in \mathbb{C}$, and both the series converges absolutely and uniformly. Then $A_j$ $=$ $B_j$ $=$ $0$, $\forall$$i$ $\in$ $\mathbb{N}$. 
\end{lemma} 
In this article $\mathbb{N}$ denotes natural numbers including $0$. \\
We did not find a proof for this result in the available literature. Therefore for completeness of presentation this lemma is proved in Section 7. Lemma 1 \ implies uniqueness of the Karp expansion (\ref{eq:3.2}) for a radiation solution. \\
In addition to expansion (\ref{eq:3.2}), we also have the following asymptotic expansion:
\begin{equation}
\label{eq:3.3}  
\psi_1(x)\sim \sqrt{\frac{2}{\pi \kappa |x|}}e^{i(\kappa |x| - \frac{\pi}{4})}\sum\limits_{j=0}^{\infty}\frac{f_j(\phi)}{|x|^{j}} \ \ \mbox{\rm for}\ \ x\in \RR^2, \ \ |x| \to \infty.
\end{equation}
However, the series in (\ref{eq:3.3}) diverges in general and in particular for $\psi_1(x) = H_0(\kappa |x|)$ (see \cite{K} for more details). \\
Theorem III \ in \cite{K} gives us the following relations:
\begin{equation}
\label{eq:3.4}  
F_0(\phi)= \frac{1}{2}[f_0(\phi) + f_0(\phi + \pi)], \;\;\;\;\;\;\;
-iG_0(\phi)=\frac{1}{2}[f_0(\phi) - f_0(\phi + \pi)] ,
\end{equation}
where $F_0, G_0$, and $f_0$ are the leading coefficients in (\ref{eq:3.2}) and (\ref{eq:3.3}) respectively.\\
In the present work we also use the following properties of $H_0$ and $H_1$:
\begin{equation}
\label{eq:3.5}  
H_\nu(r) \sim \sqrt{\frac{2}{\pi r}}e^{i(r - \frac{1}{2}\nu\pi - \frac{\pi}{4})}\sum\limits_{j=0}^{\infty}\frac{C_{\nu j}}{r^j}, \ \ as \  |x| \to \infty,
\end{equation}
\begin{equation}
C_{\nu j}=\frac{(\frac{1}{2}-\nu)_j\Gamma(\nu + j + \frac{1}{2})}{j!\Gamma(\nu + \frac{1}{2})(2i)^j}, \nonumber
\end{equation}
\begin{equation}
\label{eq:3.6}  
H_0'= -H_1, \;\;\;\;\;\;\;
H_0''=\frac{1}{r}H_1 - H_0, 
\end{equation}
\begin{equation}
\label{eq:3.7}  
\frac{H_1'}{r} = \frac{-2H_0''}{r}, \;\;\;\;\;\;\;
H_1''= -(\frac{1}{r}H_1 - H_0)',
\end{equation}
where $'$ denotes derivative with respect to $r$.\\Recall that,
\begin{align}
&\psi_1(x)=  2\int_{L}\frac{\partial G^+(x-y,\kappa)}{\partial \nu_{y}}\psi_1(y)dy,\ \ x\in V_L,  \label{eq:3.9}\\
&G^+(x,\kappa)= \frac{i}{4}H_0(\kappa|x|),\ \ x\in \RR^2, \nonumber
\end{align}
where $\psi_1$ is a radiation solution of equation (\ref{eq:1.1}), $L$ and $V_L$ are line and open half-plane in ${\cal U}$,
where $L$ is the boundary of $V_L$, $\nu$  is the outward normal to $L$ relative to $V_L$; see, for example, formula 5.84 in \cite{BR}.

\section{Determining $F_n$ and $G_n$ from $f_i$}

In this  section, we  give a method to determine $F_n(\phi)$ and $G_n(\phi)$ in (\ref{eq:3.2}), from $f_0(\phi)$,$f_0(\phi + \pi)$,...,$f_n(\phi)$,$f_n(\phi +\pi)$  values in (\ref{eq:3.3}), for $\phi$   and $(\phi + \pi)$. 

\subsection{Recursion relations for $F_n, G_n$}
\begin{lemma}
The coefficients $F_n = F_n(\phi),G_n = G_n(\phi)$ in (\ref{eq:3.2}) satisfy the following recursive relations:
\begin{equation}
\label{eq:5.1} 
F_{n+1}=\frac{-((n+1)^2G_n + G_n'')}{2\kappa(n+1)}, \;\;\;\;\;\;\;
G_{n+1}=\frac{n^2F_n + F_n''}{2\kappa(n+1)},\ \ \ n\in \NN,
\end{equation}
where $F_n'',G_n''$ denote double derivates with respect to $\phi$.
\end{lemma}
The fact that $F_n,G_n$ can be obtained recursively in terms of $F_0,G_0$ was already mentioned in \cite{K}, where formula (\ref{eq:5.1}) was given for $n=0$.\\
The derivation of relations (\ref{eq:5.1}) is as follows.

We substitute expansion (\ref{eq:3.2}) into the Helmholtz equation (\ref{eq:1.1}) (using (\ref{eq:3.8})) for $|x|>\rho$ and express the derivatives of $H_0$ and $H_1$ with respect to $r$ using the formulas (\ref{eq:3.6}) and (\ref{eq:3.7}). 

As a result, we have: 
\begin{equation}
(\Delta + \kappa^2)\psi_1(x) = \nonumber
\end{equation}
\begin{align}
\label{eq:5.7}
&H_0(\kappa r)\sum\limits_{n=0}^{\infty}\bigl(F_n'' - 2\kappa (n+1)G_{n+1} + n(n+1)F_n - nF_n \bigr)\frac{1}{r^{n+2}}\\ \nonumber
&+ H_1(\kappa r)\sum\limits_{n=0}^{\infty}\bigl(G_n'' + 2\kappa (n+1)F_{n+1} + n(n+1)G_n + (n+1)G_n \bigr)\frac{1}{r^{n+2}} = 0, 
\end{align}
\begin{center}
{$\mbox{\rm for}\ \ r > \rho.$}
\end{center}
Due to \cite{K}, both series in (\ref{eq:5.7}) converge absolutely and uniformly for $r\geqq\rho_1>\rho$. Therefore, on applying Lemma 1, we get the relations in (\ref{eq:5.1}).

\subsection{Symmetries for $F_n, G_n$}
\begin{lemma}
The coefficients $F_n, G_n$ in (\ref{eq:3.2}) have the following important symmetries:

\begin{equation}
\label{eq:5.2}  
F_n(\phi)= (-1)^nF_n(\phi + \pi), \;\;\;\;\;\;\;
G_n(\phi)=(-1)^{n+1}G_n(\phi + \pi), 
\end{equation}
\begin{center}
\begin{center}
{$\mbox{\rm for}\ \ n\in \NN.$}
\end{center}
\end{center} 
\end{lemma}
The proof of Lemma 3 \ is as follows.\\
Due to formulas (\ref{eq:3.4}), symmetries (\ref{eq:5.2}) are true for $n = 0$.
In addition, using relations (\ref{eq:5.1}) one can see that if symmetries (\ref{eq:5.2}) hold for $n$ then they are true for $(n+1)$. Hence, by induction (\ref{eq:5.2}) is true $\forall n \in \mathbb{N}$.

\begin{remark}
Suppose $\psi_1$ satisfies (\ref{eq:1.1}) and (\ref{eq:1.3}). Let $L,L^+,L^-$ be as in (\ref{eq:1.4}), $L^+,L^- \subset{\cal{U}}$ and ${0} \subset L$. Then from analyticity of $\psi_1$ in $\cal{U}$, symmetries (\ref{eq:5.2}) and expansion (\ref{eq:3.2}) it follows that $\psi_1$ on $L^+$ uniquely determines $\psi_1$ on $L^-$.
\end{remark}

\subsection{Comparison of expansions (\ref{eq:3.2}) and (\ref{eq:3.3})}

Now we equate asymptotically the series (\ref{eq:3.3}) with the convergent series (\ref{eq:3.2}) as $r \to \infty$.\\
The coefficients of $r^{-n}$ must be equal for these asymptotic expansions $\forall n \in \mathbb{N}$.

As a result, we get the following relations:
\begin{equation}
\label{eq:5.3}  
f_n(\phi)=F_n(\phi) - iG_n(\phi) + \sum\limits_{j=1}^{n}(C_{0j}F_{n-j}(\phi) - iC_{1j}G_{n-j}(\phi)),
\end{equation}
\begin{equation}
\label{eq:5.4}  
f_n(\phi + \pi)=F_n(\phi + \pi) - iG_n(\phi + \pi)  +  \sum\limits_{j=1}^{n}(C_{0j}F_{n-j}(\phi + \pi) - iC_{1j}G_{n-j}(\phi + \pi)).
\end{equation}

Proceeding from formulas (\ref{eq:5.2}), (\ref{eq:5.3}), (\ref{eq:5.4}) we obtain the following result.
\begin{theorem}
The coefficients $F_n, G_n$ in (\ref{eq:3.2}) and $f_n$ in (\ref{eq:3.3}) are related by the formulas:
\begin{equation}
\label{eq:5.5}  
F_n(\phi) - iG_n(\phi)=f_n(\phi) - \sum\limits_{j=1}^{n}(C_{0j}F_{n-j}(\phi) - iC_{1j}G_{n-j}(\phi))  ,
\end{equation}
\begin{equation}
\label{eq:5.6}  
F_n(\phi) + iG_n(\phi)=(-1)^n(f_n(\phi + \pi) -\sum\limits_{j=1}^{n}(-1)^{(n-j)}(C_{0j}F_{n-j}(\phi) + iC_{1j}G_{n-j}(\phi)) ).
\end{equation}
By these formulas, the coefficients $f_0,...,f_n$ uniquely determine the coefficients $F_0,G_0,...,F_n,G_n$ in a recursive way.
\end{theorem}
Theorem 3 \ for $n=0$ reduces to formulas (\ref{eq:3.4}).

In particular, we consider relations (\ref{eq:5.5}), (\ref{eq:5.6}) as a linear system for finding $F_n(\phi), G_n(\phi)$ from $f_n(\phi), f_n(\phi + \pi)$ and $F_0(\phi), G_0(\phi),...,F_{n-1}(\phi),G_{n-1}(\phi)$, that are found from $f_0(\phi),...,f_{n-1}(\phi)$.

\section{Determining $f_i$}
Recall that already work \cite{N4} gives formulas for finding all $f_i$ in (\ref{eq:3.3}) from $|\psi|^2 = |\psi_0 + \psi_1|^2$ on $L^+_{0,\theta}$. However, these formulas in \cite{N4} are not very simple. Therefore, in a similar way to \cite{N} (where three-dimensional case is considered) we give a very simple proof that $|\psi|^2$ on $L^+\subseteq L^+_{0, \theta}$ defined as in (\ref{eq:1.4}) where ${0}$ denotes the origin in $\RR^2$, uniquely determines  $f_j(\phi)$  in (\ref{eq:3.3}) $\forall j \in \mathbb{N}$,
under the assumption that $\theta \neq k/|k|$,  where  $\psi=\psi_0+\psi_1$,  $\psi_0$  and $\psi_1$ are solutions of equation (\ref{eq:1.1}) as in (\ref{eq:1.2}) and (\ref{eq:1.3}). \\
Let 
\begin{equation}
\label{eq:4.1}  
a(x,k)=\sqrt{|x|}(|\psi(x)|^2-1),\ \ x\in{\cal U},
\end{equation}
where $\psi=\psi_0+\psi_1$,  $\psi_0$  and $\psi_1$ are  solutions of equation (\ref{eq:1.1}) as in (\ref{eq:1.2}) and (\ref{eq:1.3}), $k$ is the wave vector in (\ref{eq:1.2}).
Then as $|x|\to +\infty$  (see \cite{N1}, \cite{N4}):
 \begin{equation}
\label{eq:4.2}  
a(x,k)=\sqrt{\frac{2}{\pi \kappa}}\bigl(e^{i(\kappa|x|-kx-\frac{\pi}{4})}f_0(\phi)+  e^{-i(\kappa|x|-kx-\frac{\pi}{4})} \overline {f_0(\phi)}\bigr)+O\left(\frac{1}{\sqrt{|x|}} \right),  
\end{equation}
uniformly in $\theta=x/|x|$;
\begin{align}
&f_0(\phi)= \frac{i\sqrt{\frac{\pi \kappa}{2}}}{D}\bigl(e^{i(ky-\kappa|y|+\frac{\pi}{4})}a(x,k)-e^{i(kx-\kappa|x|+\frac{\pi}{4})}a(y,k)\bigr)+O\left(\frac{1}{\sqrt{|x|}} \right), \label{eq:4.3}\\
&D=2sin(\tau(k\theta-\kappa)),\ \ \ \  \theta\in\mathbb S^1,\ \  \tau>0, \nonumber \\
&x,y\in L^+_{x_1, \theta},\ \ x_1=0,\ \ y=x+\tau\theta, \nonumber
\end{align}
for $D\ne 0$ for fixed $\theta$ and  $\tau$,  where $L^+_{x_1, \theta}$ is defined in  (\ref{eq:1.4}).

\begin{remark}
If an arbitrary function $a$ on  $L^+_{0, \theta}$ satisfies (\ref{eq:4.2}),  for fixed $\theta \in \mathbb{S}^1$,  $k\in\RR^2$,  $\kappa=|k|>0$,
then formula (\ref{eq:4.3}) holds.
\end{remark}

The determination of $f_0$ follows from (\ref{eq:4.3}).

Suppose that $f_0$,...,$f_n$ are determined, then the determination of $f_{n+1}$ is as follows.

Let
\begin{equation}
\label{eq:4.4}  
\psi_{1,n}(x)=\sqrt{\frac{2}{\pi \kappa |x|}}e^{i(\kappa |x| - \frac{\pi}{4})}\sum\limits_{j=0}^{n}\frac{f_j(\phi)}{|x|^{j}}, \ \ \mbox{\rm where}\ \ \theta=\frac{x}{|x|},
\end{equation}
\begin{equation}
\label{eq:4.5}  
a_n(x,k)=\sqrt{|x|}(|e^{ikx}+\psi_{1,n}(x)|^2-1),
\end{equation}
\begin{equation}
\label{eq:4.6}  
b_n(x,k)=|x|^{n+1}(a(x,k)-a_n(x,k)),
\end{equation}
where $x$ is as in (\ref{eq:4.3}), $a(x,k)$ is defined by (\ref{eq:4.1}).

We have that:
\begin{align}
&a(x,k)=\sqrt{|x|}\bigl(( e^{ikx}+\psi_{1,n}(x)+\sqrt{\frac{2}{\pi \kappa |x|}}e^{i(\kappa |x| - \frac{\pi}{4})}\frac{f_{n+1}(\phi)}{|x|^{n+1}}+ O(|x|^{-n-2-\frac{1}{2}})) \label{eq:4.7}\\
&\times(e^{-ikx}+\overline{\psi_{1,n}(x)}+\sqrt{\frac{2}{\pi \kappa |x|}}e^{-i(\kappa |x| - \frac{\pi}{4})}\frac{\overline{f_{n+1}(\phi)}}{|x|^{n+1}}+ O(|x|^{-n-2-\frac{1}{2}}))-1\bigr)   \nonumber\\
&=a_n(x,k)+  \sqrt{\frac{2}{\pi \kappa }}\bigl(e^{i(\kappa|x| - \frac{\pi}{4})}\frac{f_{n+1}(\phi)}{|x|^{n+1}}+ e^{-i(\kappa |x| - \frac{\pi}{4})}\frac{\overline{f_{n+1}(\phi)}}{|x|^{n+1}}\bigr) + O(|x|^{-n-1-\frac{1}{2}}) \nonumber
\end{align}
\begin{equation}
\label{eq:4.8}  
b_n(x,k)=\sqrt{\frac{2}{\pi \kappa }}\bigl(e^{i(\kappa |x| - \frac{\pi}{4} - kx)}{f_{n+1}(\phi)} + e^{-i(\kappa |x| - \frac{\pi}{4} - kx)}\overline{f_{n+1}(\phi)}\bigr)+O\left(\frac{1}{\sqrt{|x|}}\right),
\end{equation}
as $|x|\to +\infty$, uniformly in $\theta=x/|x|$, where $a$ is defined by (\ref{eq:4.1}).

Due to (\ref{eq:4.8}) and Remark 2, as $|x| \to +\infty$ we get:
\begin{align}
&f_{n+1}(\phi)=\frac{i\sqrt{\frac{\pi \kappa}{2}}}{D}\bigl(e^{i(ky-\kappa|y| + \frac{\pi}{4})}b_n(x,k)-e^{i(kx-\kappa|x| + \frac{\pi}{4})}b_n(y,k)\bigr)+O\left(\frac{1}{\sqrt{|x|}}\right) , \label{eq:4.9}\\
&D=2sin(\tau(k\theta-\kappa)),\ \  x,y\in L^+\subseteq L^+_{0, \theta},\ \  y=x+\tau\theta,\ \  \theta\in\mathbb S^1,\ \  \tau>0,  \nonumber
\end{align}
assuming that  $D\ne 0$ for fixed $\theta$ and  $\tau$
(where  the parameter $\tau$ can be always fixed in such a way that  $D\ne 0$, under our assumption that $\theta \neq k/|k|$).

Formulas (\ref{eq:4.1}), (\ref{eq:4.4})-(\ref{eq:4.6})  and (\ref{eq:4.9}) determine  $f_{n+1}$, give the step of induction for finding all $f_j$.

\section{Proof of Theorem 1}
\subsection{Case  $L^+\subseteq L^+_{0, \theta}$, $L^-\subseteq L^-_{0, \theta}$}
First, we give the proof for the case when $L^+\subseteq L^+_{0, \theta}$, $L^-\subseteq L^-_{0, \theta}$, where $L^+, L^+_{0, \theta}$, $L^-, L^-_{0, \theta}$ are defined in (\ref{eq:1.4}) and ${0}$ denotes the origin in $\RR^2$. \\
In this case, $|\psi|^2$ on $L^+$ and $L^-$ uniquely determines  $f_j(\phi)$ and $f_j(\phi + \pi)$ respectively in (\ref{eq:3.3}), $\forall j \in \mathbb{N}$ via formulas (\ref{eq:4.8}), (\ref{eq:4.9}),  where $\psi=\psi_0+\psi_1$; see Section 5. \\
In turn, using Theorem 3, in Section 4 we uniquely determine $F_j(\phi), F_j(\phi + \pi)$ and $G_j(\phi), G_j(\phi + \pi)$, \  $\forall j \in \mathbb{N}$. \\
Finally, the convergence of the series in (\ref{eq:3.2}) and analyticity of $\psi_1$ and  $|\psi|^2$  on~$L^+$ and $L^-$ gives us $\psi_1$ on $L^+ \cup L^-$.
\subsection{General case}
In fact, in a similar way with \cite{N}, the general case reduces to the case of Subsection 6.1 by the change of variables 
\begin{align}
&x'=x-q \label{eq:6.1}\\
&\mbox{\rm for some fixed }\ \ q\in\RR^2\ \ \mbox{\rm such that}\ \ L^+\subseteq L^+_{q, \theta},\ \  L^-\subseteq L^-_{q, \theta}.   \nonumber
\end{align}
In the new variables $x'$, we have that:
\begin{equation}
\label{eq:6.2}  
\psi_0= e^{ikq}e^{ikx'},  
\end{equation}
\begin{align}
&\psi_1\ \ \mbox{\rm   satisfies}\ \  (\ref{eq:1.3})\ \ \mbox{\rm and admits presentation}\ \  (\ref{eq:3.2}) \label{eq:6.3}\\
&\mbox{\rm with}\ \  x'\  \mbox{\rm in place of}\ \ x,\ \ \mbox{\rm with some new}\ \  r', \ F_j', \ G_j' \ \  \mbox{\rm and}\ \ f_j',   \nonumber
\end{align}
\begin{equation}
\label{eq:6.4}  
|\psi|^2=|e^{ikx'}+ e^{-ikq}\psi_1(x')|^2,  
\end{equation}
where $x'\in{\cal U'}={\cal U}-q$, and $r'$ is such that $\RR^2\setminus B_{r'} \subset {\cal U'}$;
\begin{equation}
\label{eq:6.5}  
L^+\subseteq L^+_{q, \theta}=L^+_{0, \theta} \;\;\;\;\;\;\;
L^-\subseteq L^-_{q, \theta}=L^-_{0, \theta}.   
\end{equation}
In addition,
\begin{equation}
\label{eq:6.6}  
e^{-ikq}\psi_1(x')=\sqrt{\frac{2}{\pi \kappa |x'|}}e^{i(\kappa |x'| - \frac{\pi}{4})}\sum\limits_{j=0}^{\infty}\frac{f_j''(\phi)}{|x'|^{j}}  \ \ \mbox{\rm for}\ \ x'\in \RR^2\setminus B_{r'}, \ \  
\end{equation}
where $f_j''=e^{-ikq}f_j'$. \\
We also have
\begin{equation} 
\label{eq:6.7} 
e^{-ikq}\psi_1(x')=H_0(\kappa |x'|)\sum\limits_{j=0}^{\infty}\frac{{\cal F}_j(\phi)}{|x'|^{j}} + H_1(\kappa |x'|)\sum\limits_{j=0}^{\infty}\frac{{\cal G}_j(\phi)}{|x'|^{j}},
\end{equation}
where ${\cal F}_j = e^{-ikq}F_j'$ and ${\cal G}_j = e^{-ikq}G_j'$, and the series converges absolutely and uniformly in $|x| > r'$.\\
The coefficients $F_j'$, $G_j'$ and $f_j'$ are mentioned in (\ref{eq:6.3}).

In view of (\ref{eq:6.4})-(\ref{eq:6.7}), similarly to the case of  Subsection 6.1,  it is sufficient to prove that  $|\psi|^2$ on $L^+$ and $L^-$ uniquely determines  ${\cal F}_j(\phi)$, ${\cal G}_j(\phi)$ (subsequently for $(\phi + \pi)$) in (\ref{eq:6.7}) $\forall j \in \NN$.
This determination of ${\cal F}_j$, ${\cal G}_j$ is completely similar to the determination of $F_j$ and $G_j$ in Subsection 6.1.

This completes the proof of Theorem 1.

\section{Proof of Lemma 1}
First, we recall some properties of $H_0$ and $H_1$. \\
Both $H_0$ and $H_1$ admit holomorphic continuation from $\mathbb{R}_+$ to $\mathbb{C}_- = \{z\in \mathbb{C}| \ Im(z) \\ < 0 \}$; see \cite{AS}.  \\
In view of (\ref{eq:3.6}) we also have that:
\begin{align}
\label{eq:7.1}
\frac{dH_0(z)}{dz}= -H_1(z), \ \ \ z \in \CC_-.  
\end{align}
In addition, $H_0$ has infinitely many zeroes $z_j \in \mathbb{C}_-$ such that $|z_j| \to \infty$, as $j \to \infty$; see \cite{AS} and \cite{CS}. \\ 
Using formula (\ref{eq:7.1}) we get:
\begin{align}
\label{eq:7.2}
\lim_{z\to z_j} \frac{H_0(z)}{H_1(z)}= 0 , \ \ \ \forall j\in\NN.  
\end{align}
Next, we consider: 
\begin{align}
\label{eq:7.3}
\Sigma_0(r) = \sum\limits_{j=0}^{\infty}\frac{A_j}{r^{j}},  \ \ \ \ \
\Sigma_1(r) = \sum\limits_{j=0}^{\infty}\frac{B_j}{r^{j}}, \ \ \ r \geq r_1.
\end{align}
Due to the assumptions of Lemma 1, the functions $\Sigma_0(r)$ and $\Sigma_1(r)$  admit holomorphic continuation from the interval $[r_1,+\infty)$ to the domain $\mathcal{D}_{r_1}$ where,
\begin{align}  
\mathcal{D}_\rho = \{z\in \CC \ | \  |z|\geq \rho\} \cup \{\infty\}. 
\end{align}
The rest of the proof of Lemma 1 is as follows.\\
Suppose that $\Sigma_1$ is not identically ${0}$. Then we consider the smallest $j$ such that $B_j \neq 0$, and in particular:
\begin{align}
\label{eq:7.4}
z^j\Sigma_1(z) \to B_j \ as \ z \to \infty.
\end{align}
Using formula (\ref{eq:7.4}) we get that there exists a neighbourhood of $\infty$ in $\mathcal{D}_{r_1}$ where $\infty$ is the only zero of $\Sigma_1$. \\
Then proceeding from (\ref{eq:3.8}) we get:
\begin{equation}
\label{eq:7.5}
H_0(z)\Sigma_0(z) = -H_1(z)\Sigma_1(z) \ \ \ \ in \ \ \CC_- \cap \mathcal{D}_{r_1}.  
\end{equation}
Using also (\ref{eq:7.2}) we get that:
\begin{equation}
\label{eq:7.6}
-\lim_{z\to z_j}\frac{H_0(z)\Sigma_0(z)}{H_1(z)} =  \lim_{z\to z_j}\Sigma_1(z) = 0,\ \ \ \ \forall j \in \NN, \ \ \ such \ that \ z_j \in \CC_- \cap \mathcal{D}_{r_1}.  
\end{equation}
However, formula (\ref{eq:7.6}) is in contradiction with the aforementioned conclusion that $\infty$ is the only zero of $\Sigma_1$ in $\mathcal{D}_{r_1}$. \\
Hence,
\begin{align}
\label{eq:7.7}
B_j = 0 \ \ \forall j \in \mathbb{N}, \ \ \ \ \
\Sigma_1 \equiv 0,
\end{align}
under the assumptions of Lemma 1. \\
Using (\ref{eq:3.8}), (\ref{eq:7.3}), and (\ref{eq:7.7}) we get, $A_j = 0 \ \ \forall j \in \mathbb{N}$. 

This completes the proof of Lemma 1.

\section*{Acknowledgements}

The main part of this work was fulfilled during the internship of the first author in the Centre de Math\'{e}matiques Appliqu\'{e}es of Ecole Polytechnique in May-July 2024. He also acknowledges the support provided by the Kishore Vaigyanik Protsahan Yojana fellowship.

\pagebreak

\noindent
Arjun V. Nair\\
School of Mathematics,\\
IISER Thiruvanthapuram, Thiruvananthapuram, India;\\
E-mail: arjunnair172920@iisertvm.ac.in
\\
\\
\\

\noindent
Roman G. Novikov\\
CMAP, CNRS, Ecole Polytechnique,\\
Institut Polytechnique de Paris, 91128 Palaiseau, France;\\
IEPT RAS, 117997 Moscow, Russia;\\
E-mail: novikov@cmap.polytechnique.fr

\end{document}